\documentclass[11pt,reqno]{amsart}

\usepackage[T1]{fontenc}
\usepackage[utf8]{inputenc}
\usepackage{lmodern}
\usepackage{microtype}
\usepackage[margin=1.08in]{geometry}
\usepackage{amsmath,amssymb,amsthm,mathtools}
\usepackage{xcolor}
\usepackage{graphicx}
\usepackage[colorlinks=true,linkcolor=blue!45!black,citecolor=blue!45!black,urlcolor=blue!45!black]{hyperref}

\newtheorem{theorem}{Theorem}[section]
\newtheorem{proposition}[theorem]{Proposition}
\newtheorem{lemma}[theorem]{Lemma}
\newtheorem{corollary}[theorem]{Corollary}
\theoremstyle{remark}

\newcommand{\R}{\mathbb R}
\newcommand{\Gam}{\Gamma}
\newcommand{\Lam}{\Lambda}
\newcommand{\cC}{\mathcal C}
\newcommand{\N}{\mathbb N}

\newcommand{\Sp}{S_1^+}
\newcommand{\dd}{\,\mathrm d}
\newcommand{\ud}{\mathrm d}
\newcommand{\eps}{\varepsilon}
\numberwithin{equation}{section}

\title[Odd-frequency free boundary points]{Existence of odd-frequency free boundary points \\ in the Signorini problem}
\author{Xavier Fern\'andez-Real}
\address{Institute of Mathematics, \'Ecole Polytechnique F\'ed\'erale de Lausanne, Lausanne, Switzerland}
\email{xavier.fernandez-real@epfl.ch}
\thanks{The author was supported by the Swiss National Science Foundation (SNF grant PZ00P2\_208930), by the Swiss State Secretariat for Education, Research and Innovation (SERI) under contract number MB22.00034, and by the AEI project PID2024-156429NB-I00 (Spain).}
\date{}

\begin{document}

\begin{abstract}
Up to a lower-dimensional set, the free boundary of the Signorini problem consists of points with frequency that is half-integer, even (at zero density points), or odd (at full density points). Surprisingly, despite extensive work on the structure and regularity of each of these sets, the existence of odd-frequency free boundary points has remained open. In this paper, we construct a solution with a free boundary point of frequency 3.
\end{abstract}

\maketitle

\section{Introduction}

The Signorini problem, also known as the thin (or boundary) obstacle problem, is a classical free boundary problem that goes back to the work of Signorini in linear elasticity \cite{Sig33,Sig59}: it describes the equilibrium of an elastic body resting on a rigid frictionless surface, where the region of contact is not known beforehand; see also \cite{KO88}. Its scalar version also arises in models of semipermeable membranes and boundary heat control \cite{DL76} and, as the local formulation of the obstacle problem for the half-Laplacian \cite{Sil07,CS07}, in optimal stopping and option pricing \cite{Mer76,CT04}. The regularity theory for the problem has been widely studied, with remarkable progress in recent years; see \cite{Caf79,AC04,ACS08,GP09,PSU12,KPS15,DS16,FS18,CSV20,FJ21,FR21,SY22,FT23,FS24,EY24,CC24,FY25,Kit26} and references therein, as well as the survey \cite{FR22}. 

\subsection{The Signorini problem}
We denote $X=(x',y)\in\R^n\times\R$, $B_r^+:=B_r\cap\{y>0\}$, and $B_r':=B_r\cap\{y=0\}$. We say that $u\in C(B_1)$, even in $y$, is a solution to the Signorini problem (with zero obstacle) in $B_1$ if
\begin{equation}\label{eq:signorini}
\left\{
\begin{array}{rcll}
 \Delta u&=&0&\quad\text{in }B_1^+,\\
 u&\geq&0&\quad\text{on }B_1',\\
 \partial_yu&\leq&0&\quad\text{on }B_1',\\
 u\,\partial_yu&=&0&\quad\text{on }B_1'.
\end{array}
\right.
\end{equation}
That is, $u$ is harmonic in either side of the thin space $\{y=0\}$, and it is nonnegative there. At each point of $B_1'$, either $u$ and $\partial_yu$ vanishes, but which one is unknown a priori: these are called the \emph{ambiguous boundary conditions} of Signorini. Equivalently, $u$ minimizes the Dirichlet energy $\int_{B_1}|\nabla v|^2$ among all functions $v$ with the same values on $\partial B_1$ and such that $v\geq0$ on $B_1'$; or it is the least supersolution above zero on the thin space. Solutions in $B_R$ are defined in the same way, and a global solution is a solution in $B_R$ for every $R>0$.

Solutions are $C^{1,1/2}$ in $B_1^+\cup B_1'$ (i.e., on either side of the thin space, up to $\{y = 0\}$), and this is optimal, \cite{AC04}. The thin space is divided into the \emph{contact set}, where the solution touches the obstacle, and its complement, and their interface is the \emph{free boundary},
\[
 \Lam(u):=\{x'\in B_1':u(x',0)=0\}, \qquad \Gam(u):=\partial\Lam(u),
\]
where the boundary is taken in the relative topology of $B_1'$.

\subsection{The free boundary}
Free boundary points are classified according to the rate at which the solution grows away from them, which is measured by Almgren's frequency function
\begin{equation}\label{eq:frequency}
 N(r,u):=\frac{r\int_{B_r}|\nabla u|^2}{\int_{\partial B_r}u^2}.
\end{equation}
It is nondecreasing in $r$, so that $N(0^+,u):=\lim_{r\downarrow0}N(r,u)$ is well defined. If $N(0^+,u)=\kappa$ at a contact point (say, the origin), then the blow-ups of $u$ at the origin are $\kappa$-homogeneous global solutions, and we say that the origin has \emph{frequency} $\kappa$; at free boundary points, $\kappa\in\{\tfrac32\}\cup[2,\infty)$ \cite{ACS08}. We denote by $\Gam_\kappa(u)$ the set of free boundary points of frequency $\kappa$.

In $\R^2$, that is, when $n=1$, homogeneous solutions are explicit, and their homogeneities are $2m-\frac12$, $2m$, and $2m+1$, with $m\geq1$ an integer (see, e.g., \cite[Proposition~A.1]{FS18}). By a dimension reduction argument, in any dimension these three families make up the whole free boundary, up to a lower-dimensional set:
\begin{equation}\label{eq:decomposition}
 \Gam(u)=\Gam_{\rm half}(u)\cup\Gam_{\rm even}(u)\cup\Gam_{\rm odd}(u) \cup\Gam_*(u),
\end{equation}
where
\[
 \Gam_{\rm half}(u):=\bigcup_{m\geq1}\Gam_{2m-\frac12}(u),\qquad
 \Gam_{\rm even}(u):=\bigcup_{m\geq1}\Gam_{2m}(u),\qquad
 \Gam_{\rm odd}(u):=\bigcup_{m\geq1}\Gam_{2m+1}(u),
\]
and $\Gam_*(u)$, the set of free boundary points with any other frequency, has Hausdorff dimension at most $n-2$ \cite{FS18,KW13} (it may be nonempty when $n\geq2$ \cite{FS25}). Since the free boundary is $(n-1)$-rectifiable \cite{FS18,KW13}, its $(n-1)$-dimensional part is contained in the first three sets, each of which has been extensively studied:
\begin{itemize}
\item $\Gam_{\rm half}(u)$ contains the \emph{regular points}, $\Gam_{3/2}(u)$, which form a relatively open subset of the free boundary and a smooth $(n-1)$-dimensional manifold \cite{ACS08,KPS15,DS16}. The points of frequency $\frac72$, $\frac{11}2$,\,\dots\ are less understood \cite{CSV20,SY22}.
\item $\Gam_{\rm even}(u)$ is the set of \emph{singular points}, the free boundary points at which the contact set has zero density. It is contained in a countable union of $C^1$ manifolds of dimension $n-1$ \cite{GP09}; see \cite{CSV20,FJ21,SY23} for finer results.
\item $\Gam_{\rm odd}(u)$ is the set of free boundary points at which the contact set has full density. There, the blow-up is unique, and it is a harmonic polynomial in $\{y\geq0\}$ that vanishes on the thin space \cite{FRS20}, with an algebraic rate of convergence \cite{SY23,CV25}; the set is $(n-1)$-rectifiable \cite{FS18,KW13,FRS20}.
\end{itemize}

Free boundary points in $\Gam_{\rm half}(u)$ and in $\Gam_{\rm even}(u)$ are easy to produce: already in $\R^2$, the homogeneous solutions $\operatorname{Re}\big((x_1+i|y|)^{2m-\frac12}\big)$ and $\operatorname{Re}\big((x_1+i|y|)^{2m}\big)$ have a free boundary point at the origin, of frequency $2m-\frac12$ and $2m$ respectively. Surprisingly, the situation is completely different for odd frequencies: the homogeneous solutions of degree $2m+1$ vanish identically on the thin space (in $\{y\geq0\}$ they are harmonic polynomials with nonpositive normal derivative on $\{y=0\}$, as for example $y^3/3-x_1^2y$), and therefore they do not have a free boundary. In fact, despite all the results mentioned above, not a single example of a free boundary point of odd frequency has been constructed, as discussed in \cite{GP09} and in \cite[Section~8.3]{FS18}. The existence question is posed explicitly in \cite[Question~Q2]{FR21} and is also raised in \cite{FR22,SY23,FS24,CV25}. That is, it is not even known whether $\Gam_{\rm odd}(u)$ can be nonempty.

It is therefore natural to ask whether free boundary points of odd frequency exist at all. The goal of this paper is to show that they do indeed exist.

\subsection{Main result}
We construct a free boundary point of frequency $3$, the lowest odd frequency that a free boundary point can have, in $\R^3$. From now on we take $n=2$, and we denote $X=(x',y)=(x_1,x_2,y)\in\R^3$.

\begin{theorem}\label{thm:main}
There exists a solution $u$ of \eqref{eq:signorini} such that
\[
 0\in\Gam(u)\qquad\text{and}\qquad N(0^+,u)=3.
\]
Moreover, for some $A>0$,
\begin{equation}\label{eq:main-convergence}
 \frac{u(rX)}{r^3} \to A\left(\frac{|y|^3}{3}-x_1^2|y|\right) \quad\text{locally uniformly in $\R^3$, as }r\downarrow0.
\end{equation}
\end{theorem}

\begin{figure}[t]
\centering
\includegraphics[width=0.65\textwidth]{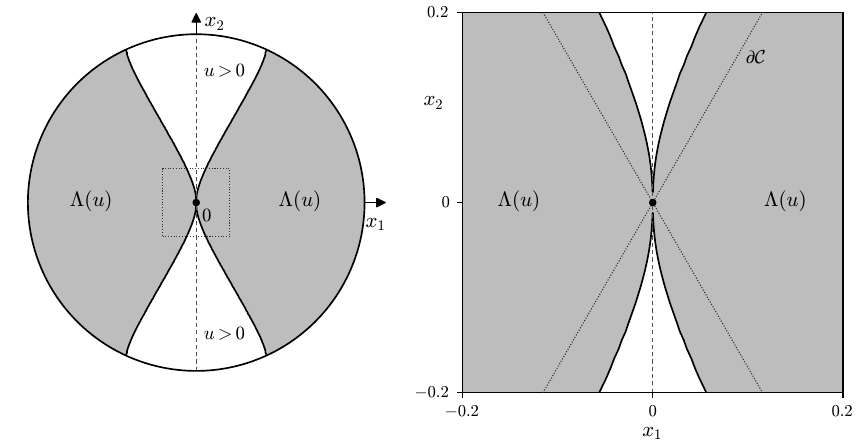}
\caption{\small The contact set $\Lam(u)$ (in grey) in $B_1'$ for the solution constructed in Section~\ref{sec:parameters}, computed numerically with $\tau=0.05$. On the right, a neighbourhood of the origin; the dotted lines are the boundary of the cone $\cC$ in \eqref{eq:cone}.}
\label{fig:contact-set}
\end{figure}

\begin{figure}[t]
\centering
\includegraphics[width=0.65\textwidth]{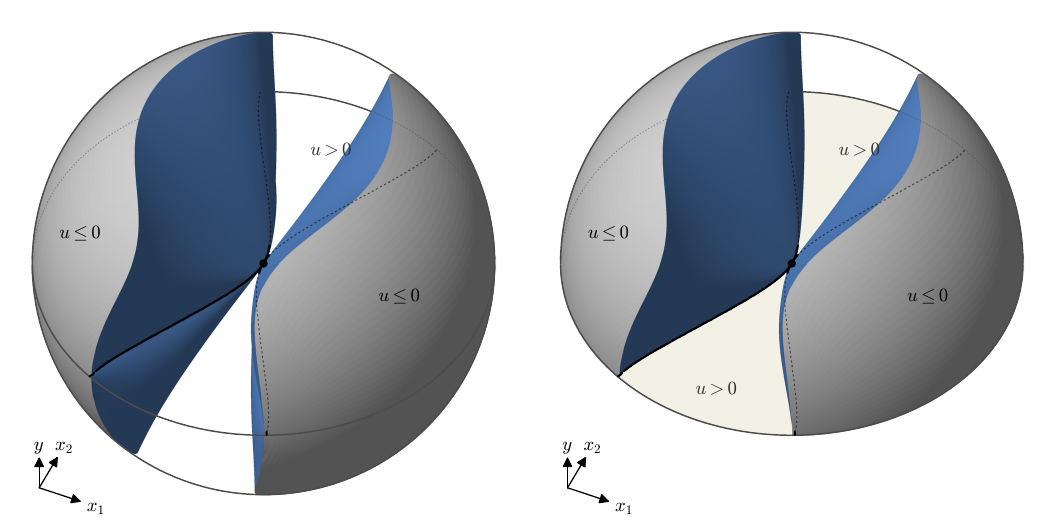}
\caption{\small The set $\{u\leq0\}$ in $B_1$ for the solution of Figure~\ref{fig:contact-set} (left), and its part in $\{y\geq0\}$ over the thin space (right). In blue, the surface $\{u=0\}$ inside the ball; in grey, the part of the boundary of the set that lies on $\partial B_1$; in black, the free boundary.}
\label{fig:negative-set}
\end{figure}

See Figure~\ref{fig:contact-set} for a depiction of the contact set of the constructed solution: the set $\{u>0\}$ consists of two regions around the spine, which meet at the origin. See also Figure~\ref{fig:negative-set} for a full space representation of the sublevel set $\{u\leq0\}$ in $B_1$ for the same solution: it consists of two regions, one on each side of $\{x_1=0\}$, which touch only at the origin.

Notice that Theorem~\ref{thm:main} also gives free boundary points of frequency $3$ in $\R^{n+1}$ for every $n\geq2$. Notice, also, that there cannot be free boundary points of odd frequency in $\R^2$; see \cite[Remark~1.2.8]{GP09} and \cite[Section~7.1]{FR22}. Indeed, in $\R^2$, homogeneous solutions with odd frequency have a strictly negative normal derivative on the thin space outside of the origin. Since rescalings converge in $C^1$ to the blow-up, a solution with such a blow-up at the origin vanishes on the thin space around it, and thus the origin is not a free boundary point.

The same observation determines the profile in \eqref{eq:main-convergence}. In $\R^3$, up to a rotation in $x'$, the $3$-homogeneous solutions are $a_1P_1+a_2P_2$ in $\{y\geq0\}$, with $a_1,a_2\geq0$ (see \cite[Lemma~5.2]{FRS20}), where
\begin{equation}\label{eq:basic-modes}
 P_1:=\frac{y^3}{3}-x_1^2y \qquad\text{and}\qquad P_2:=\frac{y^3}{3}-x_2^2y.
\end{equation}
Their normal derivative on the thin space is $-a_1x_1^2-a_2x_2^2$. If $a_1,a_2>0$, the previous argument applies, and the origin is in the interior of the contact set (see also \cite[Remarks~1.4 and~5.2]{SY23}). That is, at a free boundary point of frequency $3$ the blow-up is necessarily a multiple of $P_1$ (up to a rotation), whose normal derivative vanishes along a line: the \emph{spine} $\{x_1=0\}$.

\subsection{Idea of the proof}
Contact points of frequency $3$ in the interior of the contact set are easy to produce: the even extension of $P_1+tP_2$, plus a harmonic polynomial of higher order vanishing on $\{y=0\}$, is a solution around the origin as long as $t>0$. (See also \cite{FY25} for a classification of their counterparts at infinity: global cubic solutions with bounded positivity behave at infinity like $a_1P_1+a_2P_2+c\,|y|$ with $a_1,a_2>0$.) Our free boundary point is a limit of such points, as $t$ decreases to a given threshold. We fix $\tau>0$ small and consider solutions  in $B_1$, depending on parameters on $s$ and $t$, with prescribed boundary datum given by $P_1+tP_2-sy+\tau Q_*$ on $\partial B_1\cap\{y>0\}$, and where $Q_*$ is a suitably chosen 5-homogeneous harmonic polynomial. The two parameters, $s$ and $t$, are needed because, in the rescalings $r^{-3}u(r\,\cdot\,)$ of a solution close to $P_1$, a multiple of $y$ grows as $r\downarrow0$ and a multiple of $P_2$ is preserved. The first one is fixed by comparison: for each $t$, we take the smallest $s$ for which the origin is a contact point (Lemma~\ref{lem:critical}). Then we decrease $t$ from $8\tau$, where the solution is the polynomial itself, while the frequency at the origin is at least $3$. The danger is that the rescalings drift along the two-parameter family of profiles $a_1P_1+a_2P_2$. What prevents it is the monotonicity in $r$ of the moments $r^{-3}\int_{\partial B_1\cap\{y>0\}}u(r\theta)P_i(\theta)\dd\theta$, introduced by Figalli, Ros-Oton, and Serra in \cite[Appendix~B]{FRS20} (Lemma~\ref{lem:moment}), which keeps the solution close to $P_1$ at all scales (Proposition~\ref{prop:stability}). An argument by continuity in $t$ then shows that, at the last admissible value, the blow-up is a multiple of $P_1$ (Corollary~\ref{prop:selection}).

This does not yet give a free boundary point: the even extension of $P_1-\tau Q_*$ is a solution around the origin with blow-up $P_1$, and the origin is in the interior of its contact set. The sign of $\tau$ enters through the moment of the solution along $Q_*$, which is positive at scale one and at all smaller scales; at an interior contact point, this would force a positive normal derivative somewhere on the spine (Section~\ref{sec:fb}).

\subsection{Use of AI}\label{intro:ai}

The main idea of this manuscript is that of the author. AI was then used to refine the arguments and simplify the proofs through continuous interactions. The final version of this manuscript was then written from scratch by the author, with AI assistance to improve exposition and correct typos/mistakes. The two figures in the manuscript were numerically generated using AI.

\subsection{Organization of the paper}\label{intro:organization}
We start with some preliminaries in Section~\ref{sec:prelim}, where we also recall the monotonicity of the moments (Lemma~\ref{lem:moment}). In Section~\ref{sec:stability} we show that solutions that are close to $P_1$ and have frequency at least $3$ at the origin remain close to $P_1$ at all scales (Proposition~\ref{prop:stability}). In Section~\ref{sec:parameters} we introduce the two-parameter family of solutions and choose the parameters $s$ and $t$ (Corollary~\ref{prop:selection}). Finally, in Section~\ref{sec:fb} we show that the origin is a free boundary point, and prove Theorem~\ref{thm:main}. Appendix~\ref{app:sphere} contains two computations on the sphere that are used along the way.

\section{Preliminaries}\label{sec:prelim}

\subsection{Solutions}
We always consider solutions to the Signorini problem, \eqref{eq:signorini}, extended evenly to $y<0$. Given a solution to the Signorini problem, \eqref{eq:signorini}, in $B_R$, we will say, in addition, that it is \emph{symmetric} if it is also even in $x_1$ and $x_2$. Polynomials vanishing on $\{y=0\}$ (like $P_1$ and $P_2$ in \eqref{eq:basic-modes}) are therefore always extended evenly across $\{y=0\}$. Hence, in particular, $P_1$ and $P_2$ are global symmetric solutions.

Given a solution $u$ in $B_R$ and $r>0$, we denote
\[
 u_r(X):=\frac{u(rX)}{r^3},
\]
which is a solution in $B_{R/r}$. We further denote $\Sp:=\partial B_1\cap\{y>0\}$.

Given $g\in C(\partial B_1)$, $g\geq0$ on $\partial B_1\cap\{y=0\}$, even in $y$, there is a unique solution $u\in C(\overline{B_1})$ to the Signorini problem in $B_1$, \eqref{eq:signorini}, with $u=g$ on $\partial B_1$, which will be denoted $U[g]$. It can be obtained, for instance, as the least supersolution above the thin obstacle (see \cite[Section~3]{FR22}). By the comparison principle, we have
\begin{equation}\label{eq:comparison}
 g_1\leq g_2\ \text{ on }\partial B_1 \quad\Longrightarrow\quad U[g_1]\leq U[g_2]\ \text{ in }B_1,
\end{equation}
and therefore,
\begin{equation}\label{eq:contraction}
 \|U[g_1]-U[g_2]\|_{L^\infty(B_1)}\leq\|g_1-g_2\|_{L^\infty(\partial B_1)}.
\end{equation}
By uniqueness, if $g$ is even in $x_1$ and $x_2$, then $U[g]$ is symmetric.

Furthermore, solutions are locally $C^{1,1/2}$ in the interior of $B_1$ on either side, with estimates: there is a dimensional constant $C$ such that
\begin{equation}\label{eq:optimal-regularity}
 \|u\|_{C^{1,1/2}(\overline{B_{1/2}^+})} \leq C\|u\|_{L^\infty(B_1)}.
\end{equation}
This implies, for instance, by interpolation and \eqref{eq:optimal-regularity}:
\begin{equation}\label{eq:gradient-closeness}
 \|u-P_1\|_{L^\infty(B_\rho(X_\circ))}\leq\eps\leq1 \quad\Longrightarrow\quad \|\nabla(u-P_1)\|_{L^\infty(B_{\rho/2}^+(X_\circ))}\leq C_\rho\eps^{1/3},
\end{equation}
for any solution $u$ in $B_\rho(X_\circ)\subset B_1$ with $X_\circ\in\{y=0\}$, where $B_r^+(X_\circ):=B_r(X_\circ)\cap\{y>0\}$ and $C_\rho$ depends only on $\rho$.

Observe that, since $\partial_yP_1=-x_1^2$ on the thin space, a solution close to $P_1$ can only be positive on the thin space close to the spine $\{x_1=0\}$:

\begin{lemma}\label{lem:thin-trace}
Let $u$ be a solution in $B_1$ with $\|u-P_1\|_{L^\infty(B_1)}\leq\eps$.
\begin{enumerate}
\item There is a universal $\delta_\circ\in(0,\frac14)$ such that, if $\eps\leq\delta_\circ$, then $u(x',0)=0$ for every $x'\in B_{1/2}'$ with $|x_1|\geq\frac18$.
\item For every $\eta>0$ there is $\eps_\eta\in(0,1)$ such that, if $\eps\leq\eps_\eta$, then
\[
 \int_{B_1'}u(x',0)\dd x'\leq\eta\,\eps.
\]
\end{enumerate}
\end{lemma}

\begin{proof}
We use first \eqref{eq:gradient-closeness} with $X_\circ=0$ and $\rho=1$: given $x'\in B_{1/2}'$ with $|x_1|\geq\frac18$, for $\delta_\circ$ small enough,
\[
 \partial_yu(x',0)\leq-x_1^2+C\eps^{1/3} {}\leq-\tfrac1{64}+C\delta_\circ^{1/3}<0.
\]
Thus, $u(x',0)=0$.

(ii) It follows as in the first step, observing that $0\leq u(\cdot,0)\leq\eps$ on $B_1'$ and that the set $\{u(\cdot,0)>0\}$ can be made to have arbitrarily small measure proceeding as in (i) by using \eqref{eq:gradient-closeness}.
\end{proof}

\subsection{Almgren's frequency}
Let $u$ be a solution in $B_R$ with $u\not\equiv0$. Let us denote
\[
 H(r,u):=\int_{\partial B_r}u^2>0\qquad\text{for }r\in(0,R).
\]
We recall that the frequency function $N(r,u)$ in \eqref{eq:frequency} is nondecreasing, and we have
\begin{equation}\label{eq:H-derivative}
 \frac{\ud}{\ud r}\log\frac{H(r,u)}{r^2}=\frac{2N(r,u)}{r}.
\end{equation}
Moreover, $N(\cdot,u)\equiv\kappa$ in $(0,R)$ if and only if $u$ is $\kappa$-homogeneous in $B_R$ (see \cite{ACS08}, \cite[Section~4]{FR22}). We remark also that $N(\rho,u_r)=N(\rho r,u)$, and $N(r,\cdot)$ is continuous for the convergence in $C^1(\overline{B_r^+})$, as long as the limit does not vanish identically.

\subsection{Moments}
The following result is the main tool of the paper: it measures how a solution $u$ varies across scales against a homogeneous harmonic polynomial. It follows as in the proof of \cite[Lemma~B.1]{FRS20}:

\begin{lemma}\label{lem:moment}
Let $Q$ be a harmonic polynomial, homogeneous of degree $k$, with $Q=0$ on $\{y=0\}$. Let $u$ be a solution in $B_R$, and define
\[
 M_Q(r,u):=r^{-k}\int_{\Sp}u(r\theta)Q(\theta)\dd\theta \qquad\text{for }r\in(0,R).
\]
Then, $M_Q(\cdot,u)\in C^1((0,R))$ and
\begin{equation}\label{eq:moment-derivative}
 \frac{\ud}{\ud r}M_Q(r,u) =-r^{-2k-2}\int_{B_r'}u(x',0)\,\partial_yQ(x',0)\dd x'.
\end{equation}
\end{lemma}

\begin{proof}
Let $v(X):=u(rX)$, which is harmonic in $B_1^+$ and $C^1$ in $\overline{B_1^+}$. Since $X\cdot\nabla Q=kQ$, differentiating we get
\[
 \frac{\ud}{\ud r}M_Q(r,u) =r^{-k-1}\int_{\Sp}\big(Q\,\partial_\nu v-v\,\partial_\nu Q\big)\dd\theta,
\]
where $\nu$ is the outward unit normal to $B_1^+$. On the other hand, since $v$ and $Q$ are harmonic in $B_1^+$, Green's identity gives that the integral of $Q\,\partial_\nu v-v\,\partial_\nu Q$ on $\partial (B_1^+)$ vanishes. On the remaining boundary we have $Q=0$ and $\partial_\nu Q=-\partial_yQ$, so that
\[
 \int_{\Sp}\big(Q\,\partial_\nu v-v\,\partial_\nu Q\big)\dd\theta =-\int_{B_1'}u(rx',0)\,\partial_yQ(x',0)\dd x' =-r^{-k-1}\int_{B_r'}u(x',0)\,\partial_yQ(x',0)\dd x',
\]
where we used that $\partial_yQ$ is $(k-1)$-homogeneous.
\end{proof}

The same computation as in the previous proof applied to a $k'$-homogeneous harmonic polynomial $Q'$, with $k'\neq k$, gives
\begin{equation}\label{eq:orthogonality}
 \int_{\Sp}Q\,Q'\dd\theta=0.
\end{equation}
Observe, moreover, that for $k=3$, $M_Q(\rho,u_r)=M_Q(\rho r,u)$. In particular, since $\partial_yP_i(x',0)=-x_i^2$, Lemma~\ref{lem:moment} implies
\begin{equation}\label{eq:cubic-monotone}
 r\mapsto M_{P_i}(r,u)\quad\text{is nondecreasing, for }i=1,2.
\end{equation}

As a short consequence, we obtain the classification of cubic profiles, already contained in \cite[Lemmas~5.1 and~5.2]{FRS20} (which we include here for completeness):

\begin{lemma}\label{lem:cubic}
Let $w$ be a symmetric, $3$-homogeneous global solution. Then, $w=a_1P_1+a_2P_2$ in $\{y\geq0\}$ for some $a_1,a_2\geq0$.
\end{lemma}

\begin{proof}
By homogeneity, $M_{P_i}(r,w)$ is constant in $r$, and then \eqref{eq:moment-derivative} gives $\int_{B_1'}w(x',0)\,|x'|^2\dd x'=0$. Since $w\geq0$ on the thin space, we deduce that $w$ vanishes in $B_1'$ and, by homogeneity, on the whole thin space. Thus, the odd reflection of $w|_{\{y\geq0\}}$ is harmonic in $\R^3$ and $3$-homogeneous, that is, a harmonic cubic polynomial. Since it is odd in $y$ and even in $x_1$ and in $x_2$, it belongs to the span of $P_1$ and $P_2$, and $\partial_yw=-a_1x_1^2-a_2x_2^2\leq0$ on the thin space gives $a_1,a_2\geq0$.
\end{proof}

To detect the coefficient $a_2$ in the expansion of Lemma~\ref{lem:cubic}, we will use the moment along the orthogonal polynomial:
\begin{equation}\label{eq:Psi}
 \Psi:=P_2-\tfrac14P_1, \qquad\text{which satisfies}\qquad \int_{\Sp}P_1\Psi\dd\theta=0
\end{equation}
(see \eqref{eq:sphere-moment} in the Appendix), so that $M_\Psi(r,a_1P_1+a_2P_2)=a_2\,\|\Psi\|_{L^2(\Sp)}^2$. Since $\partial_y\Psi(x',0)=\frac14x_1^2-x_2^2$, by Lemma~\ref{lem:moment} we have
\begin{equation}\label{eq:Psi-derivative}
 \frac{\ud}{\ud r}M_\Psi(r,u) =r^{-8}\int_{B_r'}u(x',0)\big(x_2^2-\tfrac14x_1^2\big)\dd x'.
\end{equation}
Contrary to \eqref{eq:cubic-monotone}, this does not have a sign in general. However, it does have a sign when $u$ vanishes on the thin space outside of the cone
\begin{equation}\label{eq:cone}
 \cC:=\{x'\in\R^2:3x_1^2\leq x_2^2\},
\end{equation}
since
\begin{equation}\label{eq:cone-inequality}
 x_2^2-\tfrac14x_1^2\geq x_1^2\geq0 \qquad\text{for }x'\in\cC.
\end{equation}
Notice also that $4x_1^2>|x'|^2$ outside of $\cC$, and therefore
\begin{equation}\label{eq:annulus-cone}
 \big(B_{1/2}'\setminus B_{1/4}'\big)\setminus\cC \subset\big\{x'\in B_{1/2}':|x_1|>\tfrac18\big\},
\end{equation}
which is a set where the solutions in Lemma~\ref{lem:thin-trace}(i) vanish.

\section{Closeness to $P_1$ across scales}\label{sec:stability}

In general, a solution $u$ that is close to $P_1$ in $B_1$ need not be close to $P_1$ at smaller scales. In this section, we will show that this is, nonetheless, the case, when the frequency of $u$ at the origin is greater or equal than $3$. Throughout this section, the constant $\delta_\circ\in(0,\frac14)$ comes from Lemma~\ref{lem:thin-trace}: in particular,  $\delta_\circ<\|P_1\|_{L^\infty(B_1)}=\frac13$.

\begin{proposition}\label{prop:stability}
For every $\delta\in(0,\delta_\circ]$ there exists $\eps>0$ such that the following holds. Let $u$ be a symmetric solution in $B_1$ with
\[
 \|u-P_1\|_{L^\infty(B_1)}\leq\eps \qquad\text{and}\qquad N(0^+,u)\geq3.
\]
Then,
\[
 \|u_r-P_1\|_{L^\infty(B_1)}\leq\delta \qquad\text{for all }r\in(0,1].
\]
\end{proposition}

\begin{proof}
The proof follows by compactness. Assume by contradiction that for some $\delta\in(0,\delta_\circ]$ there are symmetric solutions $u_k$ in $B_1$ with $N(0^+,u_k)\geq3$ and $\|u_k-P_1\|_{L^\infty(B_1)}\leq\frac1k$ for which the conclusion fails. Since $P_1$ is $3$-homogeneous, for any $r\in(0,1)$,
\[
 \|(u_k)_r-P_1\|_{L^\infty(B_1)}=r^{-3}\|u_k-P_1\|_{L^\infty(B_r)} \leq\frac1{kr^3},
\]
and the left-hand side is continuous in $r\in(0,1]$. Hence, for $k$ large enough there is $r_k\in(0,\frac12)$, with $r_k\to0$, such that
\[
 \|(u_k)_{r_k}-P_1\|_{L^\infty(B_1)}=\delta \qquad\text{and}\qquad \|(u_k)_{r}-P_1\|_{L^\infty(B_1)}\leq\delta \quad\text{for }r\in[r_k,1].
\]
Let $w_k:=(u_k)_{r_k}$, which is a symmetric solution in $B_{1/r_k}$. Since $(w_k)_R=(u_k)_{Rr_k}$, we have
\begin{equation}\label{eq:wk-close}
 \|(w_k)_R-P_1\|_{L^\infty(B_1)}\leq\delta \quad\text{for }1\leq R\leq\frac1{r_k},
\end{equation}
with equality for $R=1$. That is, $\|w_k-P_1\|_{L^\infty(B_R)}\leq\delta R^3$. In particular, the sequence $w_k$ is locally bounded and, up to a subsequence, it converges locally in $C^1$ in $\{y\geq0\}$ to a symmetric global solution $w$ satisfying
\begin{equation}\label{eq:w-close}
 \|w-P_1\|_{L^\infty(B_R)}\leq\delta R^3\quad\text{for all }R\geq1, \qquad\text{and}\qquad \|w-P_1\|_{L^\infty(B_1)}=\delta.
\end{equation}
We will reach a contradiction by showing that $w=P_1$.

\smallskip\noindent\emph{Step 1: $w$ is a cubic profile.}
We have $w\not\equiv0$. By monotonicity, $N(\rho,w_k)=N(\rho r_k,u_k)\geq N(0^+,u_k)\geq3$, and hence $N(\rho,w)\geq3$. If we had $N(\rho_\circ,w)>3$ for some $\rho_\circ$, then \eqref{eq:H-derivative} and monotonicity of the frequency would imply that $H(R,w)$ grows at least like $R^{2+2N(\rho_\circ,w)}$ as $R\to\infty$, whereas \eqref{eq:w-close} gives $H(R,w)\leq CR^8$. Hence, $N(\cdot,w)\equiv3$, and $w$ is $3$-homogeneous, and Lemma~\ref{lem:cubic} gives
\[
 w=a_1P_1+a_2P_2\quad\text{in }\{y\geq0\},\qquad\text{with }a_1,a_2\geq0.
\]

\smallskip\noindent\emph{Step 2: almost monotone moments.}
 On the one hand, since $w_k\to w$ uniformly in $B_1$ and $w$ vanishes on the thin space,
\[
 \eta_k:=\|w_k(\cdot,0)\|_{L^\infty(B_{1/4}')}\to0.
\]
On the other hand, by \eqref{eq:wk-close} we can apply Lemma~\ref{lem:thin-trace}(i) and \eqref{eq:annulus-cone} to $(w_k)_R$ for  $1\leq R\leq\frac1{r_k}$, and in particular obtain
\begin{equation}\label{eq:fk-cone}
 w_k(\cdot,0)=0\quad\text{in }\big(B_{R_k}'\setminus B_{1/4}'\big)\setminus\cC, \qquad\text{where }R_k:=\frac1{2r_k}.
\end{equation}
Let $\rho\in[1,R_k]$. We split the integral in \eqref{eq:Psi-derivative} into $B_{1/4}'$, where we use the bound $\eta_k$, and $B_\rho'\setminus B_{1/4}'$, where by \eqref{eq:fk-cone} and \eqref{eq:cone-inequality} the integrand is nonnegative. This gives
\begin{equation}\label{eq:Psi-almost}
 \frac{\ud}{\ud\rho}M_\Psi(\rho,w_k)\geq-\eta_k\rho^{-8}.
\end{equation}
In the same way, now using $x_1^2\leq x_2^2-\frac14x_1^2$ in $\cC$,
\begin{equation}\label{eq:P1-almost}
 0\leq\frac{\ud}{\ud\rho}M_{P_1}(\rho,w_k) =\rho^{-8}\int_{B_\rho'}w_k(x',0)\,x_1^2\dd x' \leq\frac{\ud}{\ud\rho}M_\Psi(\rho,w_k)+2\eta_k\rho^{-8}.
\end{equation}

\smallskip\noindent\emph{Step 3: $a_2=0$ and $a_1=1$.}
For $Q\in\{\Psi,P_1\}$ we have $M_Q(R_k,w_k)=M_Q(\frac12,u_k)$ and $M_Q(1,w_k)\to M_Q(1,w)$. Since $u_k\to P_1$ uniformly in $B_1$, by \eqref{eq:Psi} we get
\[
 M_\Psi(R_k,w_k)\to0,\qquad M_\Psi(1,w_k)\to a_2\,\|\Psi\|^2_{L^2(\Sp)}, \qquad\text{as }k\to\infty,
\]
and
\[
 M_{P_1}(R_k,w_k)\to\|P_1\|^2_{L^2(\Sp)},\qquad M_{P_1}(1,w_k)\to a_1\,\|P_1\|^2_{L^2(\Sp)}+a_2\int_{\Sp}P_1P_2\dd\theta, \qquad\text{as }k\to\infty.
\]
Integrating \eqref{eq:Psi-almost} in $[1,R_k]$ and letting $k\to\infty$ we obtain $-a_2\,\|\Psi\|^2_{L^2(\Sp)}\geq0$, and thus $a_2=0$. Then, integrating \eqref{eq:P1-almost} in $[1,R_k]$ and letting $k\to\infty$, we get $a_1\|P_1\|^2_{L^2(\Sp)}=\|P_1\|^2_{L^2(\Sp)}$. That is, $w=P_1$, which contradicts \eqref{eq:w-close}.
\end{proof}

As a consequence of the previous result, we obtain the following. Let us denote, from now on, by $\eps_\circ$ the constant given by Proposition~\ref{prop:stability} with $\delta=\delta_\circ$, and recall the cone in \eqref{eq:cone},
\[
 \cC=\{x'\in\R^2:3x_1^2\leq x_2^2\}.
\]

\begin{corollary}\label{cor:blowup}
Let $u$ be a symmetric solution in $B_1$ with $\|u-P_1\|_{L^\infty(B_1)}\leq\eps_\circ$ and $N(0^+,u)\geq3$. Then:
\begin{enumerate}
\item $u=0$ in $B_{1/2}'\setminus\cC$.
\item $N(0^+,u)=3$, and there are $a_1>0$ and $a_2\geq0$ such that
\[
 u_r\to a_1P_1+a_2P_2 \quad\text{locally in $C^1$ in $\{y\geq0\}$,} \text{ as }r\downarrow0.
\]
\end{enumerate}
\end{corollary}

\begin{proof}
By Proposition~\ref{prop:stability}, $\|u_r-P_1\|_{L^\infty(B_1)}\leq\delta_\circ$ for all $r\in(0,1]$.

(i) Lemma~\ref{lem:thin-trace}(i) and \eqref{eq:annulus-cone} applied to $u_r$ give $u=0$ in $(B_{r/2}'\setminus B_{r/4}')\setminus\cC$ for every $r\in(0,1]$.

(ii) Since $\|u_r-P_1\|_{L^\infty(B_R)}\leq\delta_\circ R^3$ for $r\leq\frac1R$, we can take the local limit in $\R^3$ of $u_{r_j}$, for some sequence $r_j\downarrow0$, denoted $w$. Then, $w$ is a symmetric global solution with
\[
 \|w-P_1\|_{L^\infty(B_R)}\leq\delta_\circ R^3\quad\text{for all }R>0,
\]
so that $w\not\equiv0$, and for every $\rho>0$ we have $N(\rho,w)=\lim_jN(\rho r_j,u)=N(0^+,u)$. Thus, $w$ is homogeneous of degree $N(0^+,u)\geq3$, and since it grows at most cubically at infinity, the degree is $3$. By Lemma~\ref{lem:cubic}, $w=a_1P_1+a_2P_2$ with $a_1,a_2\geq0$. Moreover, $a_1>0$: otherwise $w=a_2P_2\geq0$ at the point $(\frac{\sqrt3}2,0,\frac12)$, where $P_1=-\frac13$, and this contradicts $\|w-P_1\|_{L^\infty(B_1)}\leq\delta_\circ<\frac13$.

To see that $a_1$ and $a_2$ do not depend on the sequence: by \eqref{eq:cubic-monotone}, the limits $\lim_{r\downarrow0}M_{P_i}(r,u)$ exist for $i=1,2$, and they are equal to $M_{P_i}(1,w)=a_1\int_{\Sp}P_1P_i+a_2\int_{\Sp}P_2P_i$. Since $P_1$ and $P_2$ are linearly independent, this determines $a_1$ and $a_2$.
\end{proof}

\section{The choice of the parameters}\label{sec:parameters}

Let us now introduce the family of solutions among which we will find the one from Theorem~\ref{thm:main}. Let us start by considering the harmonic polynomial
\begin{equation}\label{eq:qstar}
 Q_*:=y\,(x_1^4-12x_1^2x_2^2+8x_2^4)+y^3(2x_1^2-12x_2^2)+y^5,
\end{equation}
which is homogeneous of degree five, vanishes on $\{y=0\}$, and satisfies
\begin{equation}\label{eq:kstar}
 \partial_yQ_*(x',0)=k_*(x'):=x_1^4-12x_1^2x_2^2+8x_2^4 \leq x_1^4+8x_2^4.
\end{equation}
For the moment, we only care about \eqref{eq:kstar}; the particular choice of polynomial will be clear later. Given $\tau>0$, $s\in\R$, $t\in[-\tau,8\tau]$, let
\[
 p_{s,t}:=P_1+tP_2-sy+\tau Q_* \qquad\text{and}\qquad u_{s,t}:=U[p_{s,t}].
\]
(Recall we take the even extension in $y$.) The solutions $u_{s,t}$ are symmetric, and by the maximum principle
\begin{equation}\label{eq:above-polynomial}
 u_{s,t}\geq p_{s,t}\quad\text{in }B_1^+.
\end{equation}
Moreover, since $0\leq y\leq1$ and $|P_2|\leq y$ on $\Sp$, again by the comparison principle \eqref{eq:comparison} we have that $u_{s,t}$ is nonincreasing in $s$, and that
\begin{equation}\label{eq:comparison-t}
 u_{s+|t-t'|,t}\leq u_{s,t'}\leq u_{s-|t-t'|,t}\quad\text{in }B_1.
\end{equation}

\subsection{The parameter $s$}
For each $t$ we choose $s$ as the first value for which the solution touches the thin obstacle at the origin.

\begin{lemma}\label{lem:critical}
Let $\tau\in(0,\frac12]$ and $t\in[-\tau,8\tau]$. Then, the value
\[
 s_*(t):=\min\{s\in\R:u_{s,t}(0)=0\}
\]
is well defined, and the following holds.
\begin{enumerate}
\item $0\leq s_*(t)\leq9\tau$. Moreover, $s_*(8\tau)=0$ and $u_{0,8\tau}=p_{0,8\tau}$ in $B_1^+$.
\item $|s_*(t)-s_*(t')|\leq|t-t'|$ for all $t,t'\in[-\tau,8\tau]$.
\item The solution $u^t:=u_{s_*(t),t}$ satisfies $u^t(0)=0$ and $\nabla u^t(0)=0$.
\end{enumerate}
\end{lemma}

\begin{proof}
Let $Z$ be the closed interval given by $Z:=\{s\in\R:u_{s,t}(0)=0\}$: by monotonicity, if $s\in Z$ then $[s,\infty)\subset Z$.

If $s\in Z$, then by \eqref{eq:above-polynomial} we have that $\partial_yu_{s,t}(0)\geq\partial_yp_{s,t}(0)=-s$, and hence $s\geq0$.

On the other hand, by \eqref{eq:kstar}, on $B_1'$ we have
\[
 \partial_yp_{s,t}(x',0)=-x_1^2-tx_2^2-s+\tau k_*(x') \leq-s+\max\{-t,0\}+8\tau.
\]
Thus, for $s=9\tau$ the polynomial $p_{s,t}$ is itself a solution in $B_1$, which vanishes at the origin, and therefore $9\tau\in Z$. This shows that $s_*(t)$ is well defined, that $Z=[s_*(t),\infty)$, and that $0\leq s_*(t)\leq9\tau$. Similarly, when $t=8\tau$ and $s=0$,
\[
 \partial_yp_{0,8\tau}(x',0) \leq-x_1^2(1-\tau x_1^2)-8\tau x_2^2(1-x_2^2)\leq0 \quad\text{on }B_1',
\]
so that $u_{0,8\tau}=p_{0,8\tau}$ and $s_*(8\tau)=0$. This proves (i).

(ii) If $s\geq s_*(t)+|t-t'|$, then by \eqref{eq:comparison-t} we have $0\leq u_{s,t'}(0)\leq u_{s-|t-t'|,t}(0)=0$. Hence $s_*(t')\leq s_*(t)+|t-t'|$, and (ii) follows by exchanging $t$ and $t'$.

(iii) We have $u^t(0)=0$ by definition. The fact that $\partial_yu^t(0)=0$ follows by uniform one-sided $C^1$ regularity and the definition of $s_*(t)$.
\end{proof}

By Lemma~\ref{lem:critical}(i)-(ii) and \eqref{eq:contraction}, there is a universal constant $C_1$ such that, for all $t,t'\in[-\tau,8\tau]$,
\begin{equation}\label{eq:ut-close}
 \|u^t-P_1\|_{L^\infty(B_1)}\leq C_1\tau \qquad\text{and}\qquad \|u^t-u^{t'}\|_{L^\infty(B_1)}\leq2|t-t'|.
\end{equation}
Let us now fix a universal constant:
\begin{equation}
    \label{eq:choice-c1}
 c_1:=2^{-15}\min\big\{\|\Psi\|^2_{L^2(\Sp)},\|Q_*\|^2_{L^2(\Sp)}\big\}.
\end{equation}
From now on we assume $\tau\leq\tau_\circ$, where $\tau_\circ\in(0,\frac12]$ is a universal constant such that $C_1\tau_\circ\leq\eps_\circ$, with $\eps_\circ$ as in Corollary~\ref{cor:blowup}, and such that, by Lemma~\ref{lem:thin-trace}(ii) with $\eta=c_1/C_1$,
\begin{equation}\label{eq:thin-mass}
 \int_{B_1'}u^t(x',0)\dd x'\leq c_1\tau \qquad\text{for all }t\in[-\tau,8\tau].
\end{equation}

\subsection{The parameter $t$}
 We now choose $t$ as the smallest value such that the frequency of $u^{t'}$ at the origin is at least $3$ for all $t'\in[t,8\tau]$. Let us start with a stability property for the contact set, which will then be used in the proof of Lemma~\ref{lem:E}(iii) below: if both coefficients of the cubic blow-up are positive, then the nearby solutions that vanish at the origin must vanish on the thin space in a neighbourhood of it as well.

\begin{lemma}\label{lem:stable-contact}
Let $u$ be a symmetric solution in $B_1$ such that, for some $a_1,a_2>0$,
\[
 u_r\to a_1P_1+a_2P_2 \quad\text{locally in $C^1$ in $\{y\geq0\}$, as }r\downarrow0.
\]
Then, there exist $\rho\in(0,\frac14)$ and $\eps>0$ such that every symmetric solution $v$ in $B_1$ with
\[
 v(0)=0\qquad\text{and}\qquad\|v-u\|_{L^\infty(B_1)}\leq\eps
\]
satisfies $v=0$ in $B_\rho'$.
\end{lemma}

\begin{proof}
Let $c:=\frac14\min\{a_1,a_2\}>0$. Since $\partial_yu_r(X)=r^{-2}\,\partial_yu(rX)$, the convergence of $u_r$ gives
\[
 r^{-2}\,\partial_yu(r\theta',0)\to -a_1\theta_1^2-a_2\theta_2^2 \quad\text{uniformly in }\theta' = (\theta_1, \theta_2) \in\overline{B_1'},\text{ as }r\downarrow0.
\]
Notice that the limit is at most $-4c$ when $|\theta'|=1$, and hence there is $\rho\in(0,\frac14)$ such that
\begin{equation}\label{eq:pressure-u}
 \partial_yu(x',0)\leq-2c\,|x'|^2 \quad\text{in }B_{2\rho}'.
\end{equation}
In particular, $u=0$ in $B_{2\rho}'$, and the odd reflection of $u|_{\{y\geq0\}}$ is harmonic in $B_{2\rho}$.

Let now $v$ be a symmetric solution in $B_1$ with $v(0)=0$ and $\|v-u\|_{L^\infty(B_1)}\leq\eps$. By \eqref{eq:optimal-regularity}, $v\to u$ in $C^1(\overline{B_{1/4}^+})$ as $\eps\downarrow0$. Thus, by \eqref{eq:pressure-u}, if $\eps$ is small enough then $\partial_yv<0$ in $\overline{B_\rho'}\setminus B_{\rho/2}'$, and hence
\begin{equation}\label{eq:f-support}
 v=0\quad\text{in }\overline{B_\rho'}\setminus B_{\rho/2}'.
\end{equation}
Let $h$ be harmonic in $B_\rho^+$ with $h=v$ on $\partial B_\rho\cap\{y>0\}$ and $h=0$ on $B_\rho'$. We claim that
\begin{equation}\label{eq:pressure-h}
 \partial_yh\leq0\quad\text{in }B_\rho'.
\end{equation}
In particular, the even extension of $h$ solves the thin obstacle problem in $B_\rho$ and thus coincides with $v$, reaching the desired result.

To prove \eqref{eq:pressure-h}, notice that $v\geq h$ in $B_\rho^+$ by the comparison principle. In particular, \eqref{eq:pressure-h} holds in $B_\rho'\setminus B_{\rho/2}'$, and since $v(0)=0$ it also gives
\begin{equation}\label{eq:h-origin}
 \partial_yh(0)\leq0.
\end{equation}
On the other hand, the odd reflections of $h$ and of $u|_{\{y\geq0\}}$ are harmonic in $B_\rho$,  with $|h-u|\leq\eps$ on $\partial B_\rho$. By interior estimates for harmonic functions, denoting $\phi:=\partial_y(h-u)(\cdot,0)$ we have $|D^2\phi|\leq C\rho^{-3}\eps$ in $B_{\rho/2}'$. Since the function $\phi$ is even in $x_1$ and in $x_2$, we also have $\nabla\phi(0)=0$, and $\phi(0)=\partial_yh(0)$ because $\partial_yu(0)=0$. Thus, by \eqref{eq:pressure-u} and \eqref{eq:h-origin}, for $x'\in B_{\rho/2}'$,
\[
 \partial_yh(x',0)=\partial_yu(x',0)+\phi(x') \leq -2c|x'|^2+\phi(0) + C\rho^{-3}\eps |x'|^2 \le \partial_yh(0)-\big(2c-C\rho^{-3}\eps\big)|x'|^2\leq0
\]
if $\eps$ is small enough. This proves \eqref{eq:pressure-h}.
\end{proof}

Let us consider the set of parameters for which the frequency at the origin is at least $3$,
\[
 E:=\big\{t\in[-\tau,8\tau]:N(0^+,u^t)\geq3\big\}.
\]
By \eqref{eq:ut-close} and the choice of $\tau_\circ$, for every $t\in E$ we can apply Corollary~\ref{cor:blowup} to $u^t$. We denote by $a_1(t)>0$ and $a_2(t)\geq0$ the coefficients of its blow-up, that is, $(u^t)_r\to a_1(t)P_1+a_2(t)P_2$ as $r\downarrow0$.

\begin{lemma}\label{lem:E}
The set $E$ satisfies the following.
\begin{enumerate}
\item $E$ is closed, and $8\tau\in E$.
\item $-\tau\notin E$.
\item If $t\in E$ and $a_2(t)>0$, then there is $\sigma>0$ such that $(t-\sigma,t+\sigma)\cap[-\tau,8\tau]\subset E$.
\end{enumerate}
\end{lemma}

\begin{proof}
(i) By Lemma~\ref{lem:critical}(i), $u^{8\tau}=P_1+8\tau P_2+\tau Q_*$ in $B_1^+$, and thus $N(0^+,u^{8\tau})=3$. Moreover, $E$ is closed: by \eqref{eq:ut-close} and \eqref{eq:optimal-regularity} the function $t\mapsto N(r,u^t)$ is continuous for every $r\in(0,\frac12]$, and therefore $t\mapsto N(0^+,u^t)=\inf_{r\in(0,1/2]}N(r,u^t)$ is upper semicontinuous.

(ii) Let $t\in E$ and $s:=s_*(t)$. By Corollary~\ref{cor:blowup}(i) and \eqref{eq:Psi-derivative}--\eqref{eq:cone-inequality}, we have that $M_\Psi(\cdot,u^t)$ is nondecreasing in $(0,\frac12]$, and by Corollary~\ref{cor:blowup}(ii) it converges to $a_2(t)\|\Psi\|^2_{L^2(\Sp)}\geq0$ as $r\downarrow0$, and so $M_\Psi(\frac12,u^t)\geq0$. On the other hand, by \eqref{eq:orthogonality} and \eqref{eq:Psi}
\[
 \lim_{r\uparrow1}M_\Psi(r,u^t)=\int_{\Sp}p_{s,t}\Psi\dd\theta =t\,\|\Psi\|^2_{L^2(\Sp)},
\]
while by \eqref{eq:Psi-derivative} and \eqref{eq:thin-mass}, for $r\in[\frac12,1)$,
\[
 \Big|\frac{\ud}{\ud r}M_\Psi(r,u^t)\Big| \leq2^8\int_{B_1'}u^t(x',0)\dd x'\leq2^8c_1\tau.
\]
Therefore, $0\leq M_\Psi(\frac12,u^t)\leq t\,\|\Psi\|^2_{L^2(\Sp)}+2^7c_1\tau$, and by the choice of $c_1$ in \eqref{eq:choice-c1} we get $t>-\tau$.

(iii) Let $t\in E$ with $a_2(t)>0$. By Corollary~\ref{cor:blowup}(ii) we can apply Lemma~\ref{lem:stable-contact} to $u^t$ to obtain some $\rho$ and $\eps$, and let $t'\in[-\tau,8\tau]$ with $|t-t'|\leq\frac\eps2$. Since $u^{t'}(0)=0$ (by Lemma~\ref{lem:critical}(iii)) and $\|u^{t'}-u^t\|_{L^\infty(B_1)}\leq\eps$ (by \eqref{eq:ut-close}), we get $u^{t'}=0$ in $B_\rho'$. Then, the odd reflection of $u^{t'}|_{\{y\geq0\}}$ is harmonic in $B_\rho$. Since it is odd in $y$ and even in $x_1$ and in $x_2$, its Taylor expansion at the origin has only terms of odd degree, and the linear one vanishes by Lemma~\ref{lem:critical}(iii). Thus $|u^{t'}(X)|\leq C|X|^3$ around the origin, and hence $N(0^+,u^{t'})\geq3$. That is, $t'\in E$, and (iii) holds with $\sigma=\frac\eps2$.
\end{proof}

Let us now make the choice for $t$:

\begin{corollary}\label{prop:selection}
Let $\tau\in(0,\tau_\circ]$, and let
\[
 t_*:=\min\{t\in[-\tau,8\tau]:[t,8\tau]\subset E\}.
\]
Then, the solution $u:=u^{t_*}$ satisfies $N(0^+,u)=3$ and
\[
 u_r\to A\,P_1 \quad\text{locally in $C^1$ in $\{y\geq0\}$,} \text{ as }r\downarrow0,
\]
for some $A>0$.
\end{corollary}

\begin{proof}
The result now follows from Lemma~\ref{lem:E}, using the minimality of $t_*$, and Corollary~\ref{cor:blowup}(ii) with $A=a_1(t_*)$.
\end{proof}

\section{The origin is a free boundary point}\label{sec:fb}

Let $u=u^{t_*}$ be the solution given by Corollary~\ref{prop:selection}. To show Theorem~\ref{thm:main} it remains to see that the origin, which is a contact point, actually belongs to the free boundary. As already explained in the introduction, this cannot follow from the blow-up alone, and it is here that we use the term $\tau Q_*$ in the boundary datum, with $\tau>0$.

Let us first explain the choice of $Q_*$ in \eqref{eq:qstar}. Observe that a $5$-homogeneous harmonic polynomial vanishing on $\{y=0\}$, even in $x_1$ and in $x_2$, is determined by its normal derivative, which has the form $k(x')=\alpha x_1^4+\beta x_1^2x_2^2+\gamma x_2^4$ on $\{y=0\}$. It is:
\begin{equation}\label{eq:H-extension}
 H[k]:=y\,k-\frac{y^3}{6}\Delta k+\frac{y^5}{120}\Delta^2k.
\end{equation}
With this notation, $Q_*=H[k_*]$, where $k_*(x')=x_1^4-12x_1^2x_2^2+8x_2^4$ was defined in \eqref{eq:kstar}. If the origin were an interior contact point of $u$, then the term of order five of $u$ at the origin would be of this form, and the constraint $\partial_yu\leq0$ on the spine would only see the coefficient $\gamma$, because $k(0,x_2)=\gamma x_2^4$. The polynomial $Q_*$ is the one whose pairing on $\Sp$ with $H[k]$ gives exactly $\gamma$ (up to a universal multiplicative constant):

\begin{lemma}\label{lem:pairing}
For every $k(x')=\alpha x_1^4+\beta x_1^2x_2^2+\gamma x_2^4$ we have
\[
 \int_{\Sp}H[k]\,Q_*\dd\theta=c_*\gamma, \qquad\text{where }c_*:=\frac{16\pi}{165}.
\]
\end{lemma}

We prove Lemma~\ref{lem:pairing} in the Appendix. As before, $k_*$ changes sign, and thus the integral of a solution against $Q_*$ is not monotone in general, unless the solution vanishes on the thin space outside of the cone $\cC$ in \eqref{eq:cone}, since
\begin{equation}\label{eq:kstar-cone}
 k_*(x')=x_1^4+4x_2^2\big(2x_2^2-3x_1^2\big)\geq0 \qquad\text{for }x'\in\cC.
\end{equation}

\begin{lemma}\label{lem:moment-positive}
Let $u=u^{t_*}$ be given by Corollary~\ref{prop:selection}. Then,
\[
 M_{Q_*}(r,u)\geq\frac\tau2\,\|Q_*\|^2_{L^2(\Sp)}>0 \qquad\text{for all }r\in(0,\tfrac12].
\]
\end{lemma}

\begin{proof}
By Lemma~\ref{lem:moment},
\[
 \frac{\ud}{\ud r}M_{Q_*}(r,u)=-r^{-12}\int_{B_r'}u(x',0)\,k_*(x')\dd x'.
\]
If $r\leq\frac12$,  Corollary~\ref{cor:blowup}(i) and \eqref{eq:kstar-cone} imply that the right-hand side is nonpositive, so that $M_{Q_*}(r,u)\geq M_{Q_*}(\frac12,u)$. If  instead $r\in[\frac12,1)$, since $|k_*|\leq8$ in $B_1'$, by \eqref{eq:thin-mass} the absolute value of the right-hand side is bounded by $2^{15}c_1\tau$. Finally, \eqref{eq:orthogonality} imples
\[
 \lim_{r\uparrow1}M_{Q_*}(r,u)=\int_{\Sp}p_{s_*(t_*),t_*}\,Q_*\dd\theta =\tau\,\|Q_*\|^2_{L^2(\Sp)},
\]
and  $M_{Q_*}(\frac12,u)\geq\tau\|Q_*\|^2_{L^2(\Sp)}-2^{14}c_1\tau$. The choice of $c_1$, \eqref{eq:choice-c1}, gives the desired result.
\end{proof}

We can now give the proof of Theorem~\ref{thm:main}:

\begin{proof}[Proof of Theorem~\ref{thm:main}]
Let $\tau\in(0,\tau_\circ]$ and let $u=u^{t_*}$ be given by Corollary~\ref{prop:selection}, so that $u(0) = 0$,  $N(0^+,u)=3$, and \eqref{eq:main-convergence} holds.

Assume by contradiction that $0\notin\Gam(u)$. Then $u=0$ in $B_\rho'$ for some $\rho>0$, and the odd reflection of $u|_{\{y\geq0\}}$ is harmonic in $B_\rho$. By symmetry its Taylor expansion at the origin has only terms of odd degree, the first of which is $A\,P_1$ by Corollary~\ref{prop:selection}. Thus, for some $k(x')=\alpha x_1^4+\beta x_1^2x_2^2+\gamma x_2^4$,
\[
 u=A\,P_1+H[k]+O(|X|^7)\quad\text{in }B_\rho^+.
\]
By \eqref{eq:orthogonality} and Lemma~\ref{lem:pairing},
\[
 M_{Q_*}(r,u)=c_*\gamma+O(r^2)\quad\text{as }r\downarrow0,
\]
and then Lemma~\ref{lem:moment-positive} gives $\gamma>0$. But then
\[
 \partial_yu(0,x_2,0)=\gamma x_2^4+O(x_2^6)>0
\]
for $x_2\neq0$ small, which contradicts \eqref{eq:signorini}. Hence $0\in\Gam(u)$.
\end{proof}

\appendix

\section{Two computations}
\label{app:sphere}

For $\alpha_1,\alpha_2,\alpha_3\in\N_{\geq0}$ we have
\begin{equation}\label{eq:sphere-moment}
 \int_{\Sp}x_1^{2\alpha_1}x_2^{2\alpha_2}y^{2\alpha_3}\dd\theta =2\pi\,\frac{(2\alpha_1-1)!!\,(2\alpha_2-1)!!\,(2\alpha_3-1)!!} {(2\alpha_1+2\alpha_2+2\alpha_3+1)!!},
\end{equation}
with the convention $(-1)!!=1$.

\begin{proof}[Proof of \eqref{eq:Psi}]
Expanding $P_1^2$ and $P_1P_2$ and using \eqref{eq:sphere-moment},
\[
 \int_{\Sp}P_1^2\dd\theta =2\pi\Big(\frac{1}{63}-\frac{2}{105}+\frac{1}{35}\Big)=\frac{16\pi}{315} \qquad\text{and}\qquad
 \int_{\Sp}P_1P_2\dd\theta =2\pi\Big(\frac{1}{63}-\frac{2}{105}+\frac{1}{105}\Big)=\frac{4\pi}{315},
\]
so that $\int_{\Sp}P_1(P_2-\frac14P_1)\dd\theta=0$.
\end{proof}

\begin{proof}[Proof of Lemma~\ref{lem:pairing}]
By linearity, it is enough to consider the three cases
\[
 H[x_1^4]=yx_1^4-2y^3x_1^2+\tfrac15y^5,\quad
 H[x_1^2x_2^2]=yx_1^2x_2^2-\tfrac13y^3(x_1^2+x_2^2)+\tfrac1{15}y^5,\quad
 H[x_2^4]=yx_2^4-2y^3x_2^2+\tfrac15y^5.
\]
We multiply each of them by $Q_*$ in \eqref{eq:qstar} and integrate by means of \eqref{eq:sphere-moment}. Grouping the result according to the three terms of $H[k]$, and denoting $\omega:=2\pi/11!!$, we find
\begin{align*}
 \int_{\Sp}H[x_1^4]\,Q_*\dd\theta&=\omega\,(24-48+24)=0,\\
 \int_{\Sp}H[x_1^2x_2^2]\,Q_*\dd\theta&=\omega\,(-48+40+8)=0,\\
 \int_{\Sp}H[x_2^4]\,Q_*\dd\theta&=\omega\,(192+288+24) =\frac{16\pi}{165}.
\end{align*}
The lemma follows.
\end{proof}


\begin{thebibliography}{FRS20}

\bibitem[AC04]{AC04}
I. Athanasopoulos, L. A. Caffarelli, \emph{Optimal regularity of lower-dimensional obstacle problems}, Zap. Nauchn. Sem. POMI \textbf{310} (2004), 49--66; English translation in J. Math. Sci. (N.Y.) \textbf{132} (2006), 274--284.

\bibitem[ACS08]{ACS08}
I. Athanasopoulos, L. A. Caffarelli, S. Salsa, \emph{The structure of the free boundary for lower dimensional obstacle problems}, Amer. J. Math. \textbf{130} (2008), 485--498.

\bibitem[Caf79]{Caf79}
L. A. Caffarelli, \emph{Further regularity for the Signorini problem}, Comm. Partial Differential Equations \textbf{4} (1979), 1067--1075.

\bibitem[CS07]{CS07}
L. A. Caffarelli, L. Silvestre, \emph{An extension problem related to the fractional Laplacian}, Comm. Partial Differential Equations \textbf{32} (2007), 1245--1260.

\bibitem[CC24]{CC24}
M. Carducci, R. Colombo, \emph{Generic regularity of free boundaries in the obstacle problem for the fractional Laplacian}, preprint arXiv:2412.16066 (2024).

\bibitem[CV25]{CV25}
M. Carducci, B. Velichkov, \emph{An epiperimetric inequality for odd frequencies in the thin obstacle problem}, J. Funct. Anal. \textbf{289} (2025), 111115.

\bibitem[CSV20]{CSV20}
M. Colombo, L. Spolaor, B. Velichkov, \emph{Direct epiperimetric inequalities for the thin obstacle problem and applications}, Comm. Pure Appl. Math. \textbf{73} (2020), 384--420.

\bibitem[CT04]{CT04}
R. Cont, P. Tankov, \emph{Financial Modelling with Jump Processes}, Chapman \& Hall/CRC Financial Mathematics Series, Chapman \& Hall/CRC, Boca Raton, FL, 2004.

\bibitem[DS16]{DS16}
D. De Silva, O. Savin, \emph{Boundary Harnack estimates in slit domains and applications to thin free boundary problems}, Rev. Mat. Iberoam. \textbf{32} (2016), 891--912.

\bibitem[DL76]{DL76}
G. Duvaut, J.-L. Lions, \emph{Inequalities in Mechanics and Physics}, Springer, Berlin, 1976.

\bibitem[EY24]{EY24}
S. Eberle, H. Yu, \emph{Compact contact sets of sub-quadratic solutions to the thin obstacle problem}, Adv. Math. \textbf{444} (2024), 109635.

\bibitem[Fer22]{FR22}
X. Fern\'andez-Real, \emph{The thin obstacle problem: a survey}, Publ. Mat. \textbf{66} (2022), 3--55.

\bibitem[FJ21]{FJ21}
X. Fern\'andez-Real, Y. Jhaveri, \emph{On the singular set in the thin obstacle problem: higher order blow-ups and the very thin obstacle problem}, Anal. PDE \textbf{14} (2021), 1599--1669.

\bibitem[FR21]{FR21}
X. Fern\'andez-Real, X. Ros-Oton, \emph{Free boundary regularity for almost every solution to the Signorini problem}, Arch. Ration. Mech. Anal. \textbf{240} (2021), 419--466.

\bibitem[FT23]{FT23}
X. Fern\'andez-Real, C. Torres-Latorre, \emph{Generic regularity of free boundaries for the thin obstacle problem}, Adv. Math. \textbf{433} (2023), 109323.

\bibitem[FY25]{FY25}
X. Fern\'andez-Real, H. Yu, \emph{Global solutions to the thin obstacle problem with superquadratic growth}, Math. Ann., to appear.

\bibitem[FRS20]{FRS20}
A. Figalli, X. Ros-Oton, J. Serra, \emph{Generic regularity of free boundaries for the obstacle problem}, Publ. Math. Inst. Hautes \'Etudes Sci. \textbf{132} (2020), 181--292.

\bibitem[FS18]{FS18}
M. Focardi, E. Spadaro, \emph{On the measure and the structure of the free boundary of the lower dimensional obstacle problem}, Arch. Ration. Mech. Anal. \textbf{230} (2018), 125--184.

\bibitem[FS25]{FS25}
F. Franceschini, O. Savin, \emph{Solutions to the thin obstacle problem with non-2D frequency}, J. Eur. Math. Soc., to appear.

\bibitem[FS24]{FS24}
F. Franceschini, J. Serra, \emph{Free boundary partial regularity in the thin obstacle problem}, Comm. Pure Appl. Math. \textbf{77} (2024), 630--669.

\bibitem[GP09]{GP09}
N. Garofalo, A. Petrosyan, \emph{Some new monotonicity formulas and the singular set in the lower dimensional obstacle problem}, Invent. Math. \textbf{177} (2009), 415--461.

\bibitem[KO88]{KO88}
N. Kikuchi, J. T. Oden, \emph{Contact Problems in Elasticity: A Study of Variational Inequalities and Finite Element Methods}, SIAM Studies in Applied Mathematics, vol.~8, SIAM, Philadelphia, 1988.

\bibitem[Kit26]{Kit26}
S. Kitano, \emph{Ellipsoidal positivity sets for fractional obstacle problems with quadratic forcing}, preprint arXiv:2609.00703 (2026).

\bibitem[KPS15]{KPS15}
H. Koch, A. Petrosyan, W. Shi, \emph{Higher regularity of the free boundary in the elliptic Signorini problem}, Nonlinear Anal. \textbf{126} (2015), 3--44.

\bibitem[KW13]{KW13}
B. Krummel, N. Wickramasekera, \emph{Fine properties of branch point singularities: two-valued harmonic functions}, preprint arXiv:1311.0923 (2013).

\bibitem[Mer76]{Mer76}
R. C. Merton, \emph{Option pricing when underlying stock returns are discontinuous}, J. Financ. Econ. \textbf{3} (1976), 125--144.

\bibitem[PSU12]{PSU12}
A. Petrosyan, H. Shahgholian, N. Uraltseva, \emph{Regularity of Free Boundaries in Obstacle-Type Problems}, Graduate Studies in Mathematics, vol.~136, American Mathematical Society, Providence, RI, 2012.

\bibitem[SY22]{SY22}
O. Savin, H. Yu, \emph{Half-space solutions with $7/2$ frequency in the thin obstacle problem}, Arch. Ration. Mech. Anal. \textbf{246} (2022), 397--474.

\bibitem[SY23]{SY23}
O. Savin, H. Yu, \emph{Contact points with integer frequencies in the thin obstacle problem}, Comm. Pure Appl. Math. \textbf{76} (2023), 4048--4074.

\bibitem[Sig33]{Sig33}
A. Signorini, \emph{Sopra alcune questioni di elastostatica}, Atti Soc. It. Progr. Sc. \textbf{21} (1933), 143--148.

\bibitem[Sig59]{Sig59}
A. Signorini, \emph{Questioni di elasticit\`a non linearizzata e semilinearizzata}, Rend. Mat. e Appl. \textbf{18} (1959), 95--139.

\bibitem[Sil07]{Sil07}
L. Silvestre, \emph{Regularity of the obstacle problem for a fractional power of the Laplace operator}, Comm. Pure Appl. Math. \textbf{60} (2007), 67--112.

\end{thebibliography}
\end{document}